\documentclass[fullpage,12pt]{article}
\usepackage{graphicx,graphics,latexsym,verbatim,epsfig}

\newtheorem{theorem}{Theorem}
\newtheorem{corollary}{Corollary}

\begin{document}
\title{Diagonal Stability for a Class of Interconnected Passive Systems}
\author{ Murat Arcak\\
Department of Electrical, Computer, and Systems Engineering\\
    Rensselaer Polytechnic Institute \\
Troy, NY 12180  \\
Email: {\tt arcakm@rpi.edu}}
\date{\today } \maketitle

\begin{abstract}
We consider a class of matrices with a specific structure that
arises, among other examples, in dynamic models for biological
regulation of enzyme synthesis \cite{TysOth78}. We first show that a
stability condition given in \cite{TysOth78} is in fact a necessary
and sufficient condition for diagonal stability of this class of
matrices. We then revisit a recent generalization of \cite{TysOth78}
to nonlinear systems given in \cite{Son05}, and recover the same
stability condition using our diagonal stability result. Unlike the
input-output based arguments employed in \cite{Son05}, our proof
gives a procedure to construct a Lyapunov function. Finally we study
static nonlinearities that appear in the feedback path, and give a
stability condition  that mimics the Popov criterion.
\end{abstract}

\subsection*{Main Result}
The results of this note were triggered by the recent paper
\cite{Son05} and by several discussions with its author. We give our
main diagonal stability result in Theorem \ref{main} below, and
present its implications for stability of a class of interconnected
systems in the form of corollaries to this theorem.

\begin{theorem} \label{main}
 A matrix of the form
\begin{equation}\label{Amatrix}
A=\left[\begin{array}{ccccc} -1  & 0 & \cdots & 0 & -\gamma_1 \\ \gamma_2 & -1 &\ddots & & 0 \\
0 & \gamma_3 & -1 & \ddots & \vdots \\ \vdots & \ddots & \ddots &
\ddots & 0 \\ 0 & \cdots & 0 & \gamma_n & -1 \end{array} \right]
\quad \gamma_i>0, \ i=1,\cdots,n,
\end{equation}
 is diagonally stable; that is, it
satisfies
$$
DA+A^TD<0
$$ for some diagonal matrix $D>0$,  if and only if
\begin{equation}\label{eduardo}
\gamma_1\dots\gamma_n<\sec(\pi/n)^n.\end{equation} \hfill $\Box$
\end{theorem}

The proof is given in an appendix. The results surveyed in
\cite{Red85,bhaya} for diagonal stability of various classes of
matrices do not encompass the specific structure exhibited by
(\ref{Amatrix}). In particular, the sign reversal for $\gamma_1$ in
(\ref{Amatrix}) rules out the ``M-matrix" condition, which is
applicable when all off-diagonal terms are nonnegative.
\smallskip

\begin{figure}[h]
\vspace{-2cm}
\begin{center}
\mbox{}\setlength{\unitlength}{1.2mm}
\begin{picture}(70,50)
\put(-10,0){\psfig{figure=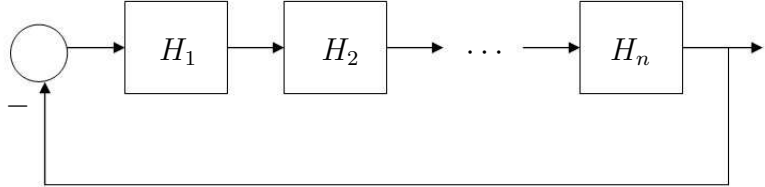,width=90\unitlength}}
\put(-7,14){$\displaystyle -$}\put(44,20){$\displaystyle \cdots$}
\put(10,20){$\displaystyle H_1$} \put(28,20){$\displaystyle H_2$}
\put(60,20){$\displaystyle H_n$}
\end{picture}
\vspace{-.5cm} \caption{Feedback interconnection for Corollary
\ref{one}.} \label{nest}
\end{center}
\end{figure}

We now apply  Theorem \ref{main}  to characterize the stability of
the feedback interconnection  in Figure \ref{nest}. When each block
$H_i$ is a first-order linear system with transfer function
$H_i(s)=\gamma_i/(\tau_is+1)$, $\gamma_i>0$, $\tau_i>0$, then a
state-space representation of the interconnection would be obtained
from the $A$-matrix in (\ref{Amatrix}) by multiplying its $i$th row
 by $1/\tau_i$ for $i=1,\cdots,n$. Since row multiplications by positive
 constants do not change diagonal stability, Theorem
 \ref{main} recovers the result in \cite{TysOth78} which states that if (\ref{eduardo}) holds then the feedback
 interconnection in Figure \ref{nest} is asymptotically stable.
 Theorem \ref{main}
 further shows that stability can be proven with a Lyapunov
 function $V=x^TPx$ in which $P$ is diagonal.
\smallskip


The linear result in \cite{TysOth78} has been extended in
\cite{Son05} to the situation where $H_i$'s are not restricted to be
linear, and instead, characterized by the {\it output feedback
 passivity} (OFP) \cite{sepulchre} (a.k.a. {\it output strict passivity} \cite{vds}) property:
\begin{equation}\label{IO}
-\beta \le -\|y_i\|^2+\gamma_i <u_i,y_i>
\end{equation}
where $\|\cdot\|$ and $<\cdot,\cdot>$ denote, respectively, the norm
and inner product in the extended $L_2$ space, and $\beta \ge 0$
represents the bias  due to initial conditions.
Using this property, \cite{Son05} proves that the {\it secant
condition} (\ref{eduardo}) insures stability of the feedback
interconnection in Figure \ref{one}.
\smallskip

Unlike the input-output proof given in \cite{Son05}, we now assume
that a {\it storage function} $V_i$ is available for each block in
Figure \ref{one}, and show that a weighted sum of these $V_i$'s,
\begin{equation}\label{lyap}V=\sum_{i=1}^{n} d_iV_i,\end{equation} where $d_i>0$ are chosen following the procedure
below,
 is a Lyapunov function for
the closed-loop system. Indeed, a storage function verifying the OFP
property (\ref{IO}) satisfies
\begin{equation}\label{ofp}
\dot{V}_i\le -y_i^2+\gamma_i u_iy_i
\end{equation}
which, when substituted in (\ref{lyap}) along with the
interconnection conditions $$u_1=-y_n, \quad u_i=y_{i-1}, \
i=2,\cdots n,$$ results in
\begin{equation}\label{RHS}
\dot{V}\le -y^TDAy
\end{equation}
where $A$ is as in (\ref{Amatrix}), and $D$ is a diagonal matrix
comprising of the coefficients $d_i$ in (\ref{lyap}). It then
follows from Theorem \ref{main} that positive $d_i$'s that render
the right-hand side of (\ref{RHS}) negative definite indeed exist if
(\ref{eduardo}) holds:

\begin{corollary} \label{one}
Consider the feedback interconnection in Figure \ref{nest} and let
$u_i$, $x_i$ and $y_i$ denote the input, state vector, and output of
each block $H_i$. Suppose, further, there exist $C^1$ storage
functions $V_i(x_i)$, satisfying (\ref{ofp}) with $\gamma_i>0$ along
the state trajectories of  each block. Under these conditions, if
 (\ref{eduardo}) holds then there exist $d_i>0$,
$i=1,\cdots,n$, such that the Lyapunov function (\ref{lyap})
satisfies
$$
\dot{V}=\sum_{i=1}^{n} d_i\dot{V_i}\le -\epsilon |y|^2
$$
for some $\epsilon>0$. \hfill $\Box$
\end{corollary}

Corollary \ref{one} still holds when some of the blocks are static
nonlinearities satisfying the sector condition
\begin{equation}\label{sectorg}
0\le -y_i^2+\gamma_i u_iy_i, \quad \gamma_i>0,
\end{equation}
rather than the dynamic property (\ref{ofp}). To see this we let
$\mathcal{I}$ denote the subset of indices $i$ which correspond to
dynamic blocks $H_i$ satisfying (\ref{ofp}), and employ the Lyapunov
function
\begin{equation}
V=\sum_{i\in \mathcal{I}} d_iV_i.
\end{equation}
For the static blocks, that is $H_i$, $i\notin \mathcal{I}$, we note
from (\ref{sectorg}) that the sum
\begin{equation}
\sum_{\begin{array}{c}i=1 \\ i\notin \mathcal{I} \end{array}}^n d_i
(-y_i^2+\gamma_i u_iy_i) \quad d_i>0
\end{equation}
is nonnegative and, hence,
\begin{equation}
\dot{V}\le \sum_{i\in \mathcal{I}}
d_i\dot{V}_i+\sum_{\begin{array}{c}i=1 \\ i\notin \mathcal{I}
\end{array}}^n d_i (-y_i^2+\gamma_i u_iy_i)\le \sum_{i=1}^n d_i(-y_i^2+\gamma_i
u_iy_i)=-y^TDAy
\end{equation}
as in (\ref{RHS}). Then, as in Corollary \ref{one}, condition
(\ref{eduardo}) insures existence of a $D>0$ such that $\dot{V}\le
-\epsilon |y|^2$ for some $\epsilon>0$.

\subsection*{A Popov Criterion}
\begin{figure}[h]
\vspace{-.5cm}
\begin{center}
\mbox{}\setlength{\unitlength}{1.2mm}
\begin{picture}(70,50)
\put(-11,10){\psfig{figure=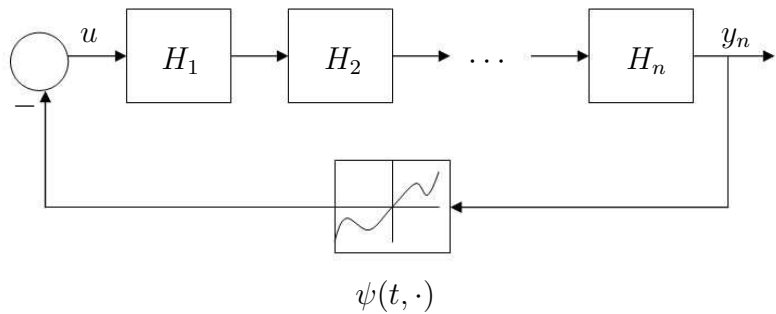,width=90\unitlength}}
\put(42.5,34){$\displaystyle \cdots$} \put(8.5,34){$\displaystyle
H_1$} \put(26.5,34){$\displaystyle H_2$} \put(60,34){$\displaystyle
H_n$} \put(-.5,37){$\displaystyle u$}\put(70.5,37){$\displaystyle
y_n$}\put(-8,29){$\displaystyle -$}\put(30,8){$\displaystyle
\psi(t,\cdot)$}
\end{picture}
\vspace{-1cm} \caption{The feedback interconnection for Corollary
\ref{four}.} \label{nest3}
\end{center}
\end{figure}

A special case of interest is the feedback interconnection
  in Figure \ref{nest3}, where $H_i$, $i=1,\cdots,n$, are dynamic
  blocks as in (\ref{ofp}), and the feedback nonlinearity
$\psi(t,\cdot)$ satisfies the {\it sector} property:
\begin{equation}\label{sect}
0\le y_n \psi(t,y_n) \le \kappa y_n^2,
\end{equation}
rewritten here as
\begin{equation}\label{sect2}
0\le -\psi(t,y_n)^2 +{\kappa}\psi(t,y_n)y_n.
\end{equation}
If we treat the feedback nonlinearity as a new block
$y_{n+1}=\psi(t,y_n)$, and note from (\ref{sect2}) that it satisfies
(\ref{sectorg}) with $\gamma_{n+1}={\kappa}$, we obtain from
Corollary \ref{one} and the ensuing discussion the stability
condition:
\begin{equation}\label{sectorcond}
\kappa\gamma_1\cdots \gamma_n <\sec(\pi/(n+1))^{(n+1)}.
\end{equation}

This condition, however, may be conservative because it does not
exploit the static nature of the feedback nonlinearity. Indeed,
using the Popov Criterion, the authors of \cite{TysOth78} obtained a
relaxed condition in which $n+1$ in the right-hand side of
(\ref{sectorcond}) is reduced to $n$ when $H_i$'s are first-order
linear blocks of the form $H_i(s)=\gamma_i/(\tau_is+1)$ and the
feedback nonlinearity is time-invariant.

\begin{figure}[h]
\begin{center}
\mbox{}\setlength{\unitlength}{1.2mm}
\begin{picture}(70,50)
\put(-11,10){\psfig{figure=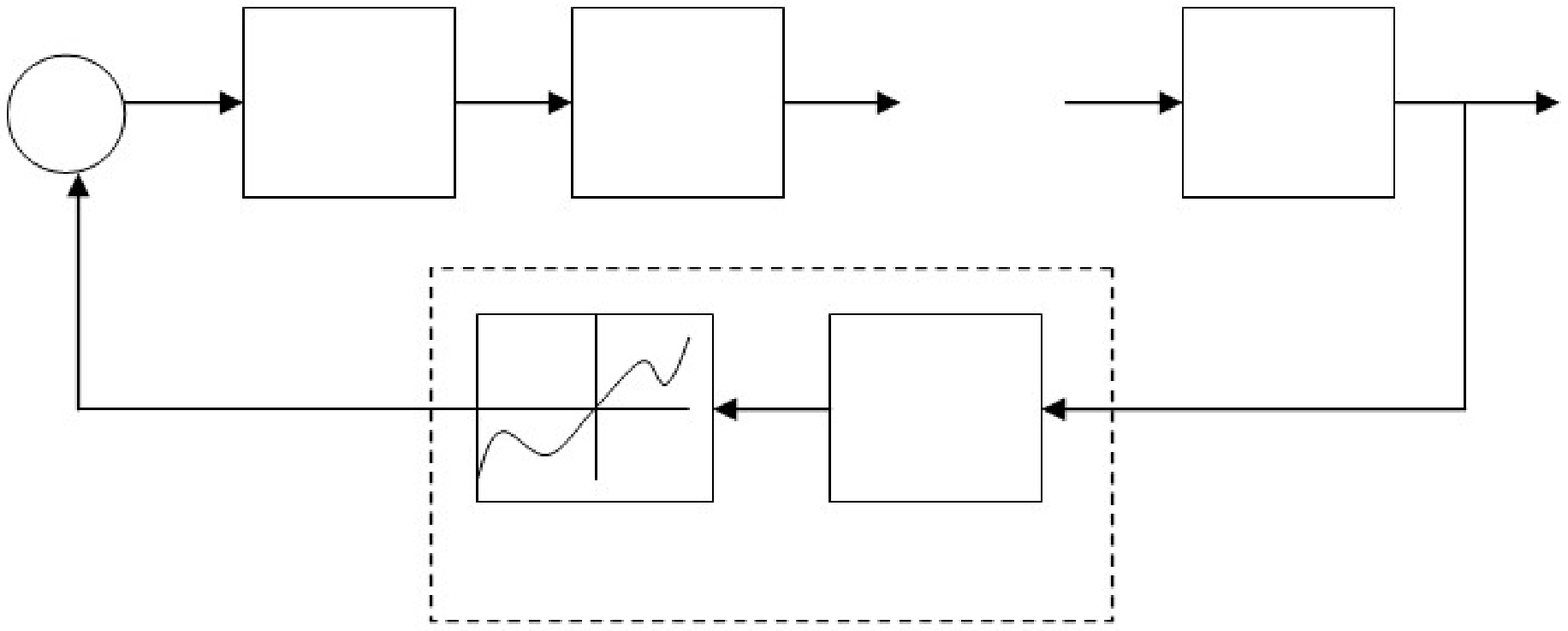,width=90\unitlength}}
\put(43.5,40){$\displaystyle \cdots$} \put(10.5,40){$\displaystyle
H_1$} \put(27.5,40){$\displaystyle H_2$} \put(59,40){$\displaystyle
H_{n-1}$} \put(2.5,43){$\displaystyle u$}\put(71,43){$\displaystyle
y_{n-1}$}\put(-8,34){$\displaystyle -$}\put(23,16){$\displaystyle
\psi(\cdot)$}\put(41.5,23.5){$\displaystyle
H_n$}\put(33,21){$\displaystyle y_n$}\put(12,21){$\displaystyle
\tilde{y}_n$}\put(33,9){$\displaystyle \tilde{H}_n$}
\end{picture}
\vspace{-1cm} \caption{An equivalent representation of the feedback
system in Figure \ref{nest3}. When $H_n$ is a linear block
$H_n(s)=\frac{\gamma_n}{\tau_ns+1}$, its series interconnection with
the $[0,\kappa]$ sector nonlinearity $\psi(\cdot)$ constitutes a
dynamic block $\tilde{H}_n$ which satisfies (\ref{ofp}) with
$\tilde{\gamma}_n=\kappa\gamma_n$.} \label{nest4}
\end{center}
\end{figure}

To extend this result to the case where $H_i$'s are OFP as in
(\ref{ofp}), and not necessarily linear,  we recall that the main
premise of the Popov Criterion is that a time-invariant sector
nonlinearity, when cascaded with a first-order, stable, linear block
preserves its passivity properties. This means that, by only
restricting $H_n$ to be linear, and combining it with the feedback
nonlinearity as in Figure \ref{nest4}, the relaxed sector condition
of \cite{TysOth78} holds even if $H_1,\cdots, H_{n-1}$ are
nonlinear:
\begin{corollary}\label{four} Consider the feedback interconnection in Figure \ref{nest3} where
 $H_i$, $i=1,\cdots, n-1$, satisfy (\ref{ofp}) with
$C^1$ storage functions $V_i$ and $\gamma_i>0$, $H_n$ is a linear
block with transfer function
\begin{equation}\label{lastblock}
H_n(s)=\frac{\gamma_n}{\tau_ns+1}, \quad \tau_n>0, \gamma_n>0,
\end{equation}
the feedback nonlinearity $\psi(\cdot)$ is time-invariant  and
satisfies the sector property (\ref{sect}). Under these assumptions,
if
\begin{equation}\label{sectorcond2}
\kappa\gamma_1\cdots \gamma_n <\sec(\pi/n)^{n},
\end{equation}
then there exists a Lyapunov function of the form
\begin{equation}
V=\sum_{i=1}^{n-1}d_iV_i+d_n \int_0^{y_n}\psi(\sigma)d\sigma, \quad
d_i>0, \ i=1,\cdots,n,
\end{equation}
satisfying
$$
\dot{V}\le -\epsilon |(y_1,\cdots, y_{n-1}, \psi(y_n))|^2
$$
for some $\epsilon>0$. \hfill $\Box$
\end{corollary}

\noindent {\bf Proof:} Rather than treat $H_n$ and $\psi(\cdot)$ as
separate blocks, we combine them as in Figure \ref{nest4}:
\begin{equation}\label{yn}
\tilde{H}_n: \ \left\{
\begin{array}{lll} \tau_n\dot{y}_n &=& -y_n+\gamma_ny_{n-1} \\ \tilde{y}_n&=&\psi(y_n),\end{array}  \right.
\end{equation}
and define
\begin{equation}
V_n=\kappa\tau_n \int_0^{y_n}\psi(\sigma)d\sigma
\end{equation}
which, from (\ref{yn}), satisfies
\begin{equation}
\dot{V}_n=-\kappa y_n\psi(y_n) +\kappa \gamma_n \psi(y_n)y_{n-1}.
\end{equation}
Because $-\kappa y_n\psi(y_n)\le -\psi(y_n)^2$ from  (\ref{sect2}),
we conclude
\begin{equation}
\dot{V}_n\le -\psi(y_n)^2 +\kappa \gamma_n
\psi(y_n)y_{n-1}=-\tilde{y}_n^2+\kappa \gamma_n \tilde{y}_ny_{n-1},
\end{equation}
which shows that $\tilde{H}_n$ is OFP as in (\ref{ofp}), with
 $\tilde{\gamma}_n=\gamma_n\kappa$. The result then follows from Corollary
\ref{one}. \hfill $\Box$
\bigskip

Corollary \ref{four} can be further generalized to the situation
where other nonlinearities exist in between the blocks $H_i$,
$i=1,\cdots,n$, in Figure \ref{nest3}. If such a nonlinearity is in
the sector $[0,\kappa_{i+1}]$, and is preceded by a linear block
$H_i(s)=\frac{\gamma_i}{\tau_is+1}$, then the two can be treated as
a single block with $\tilde{\gamma}_i=\kappa_{i+1}\gamma_i$, thus
 reducing $n$ in the right-hand side of (\ref{eduardo}).

\subsection*{The Shortage of Passivity in a Cascade of OFP Systems} 
\begin{figure}[h]
\vspace{-2.5cm}
\begin{center}
\mbox{}\setlength{\unitlength}{1.2mm}
\begin{picture}(70,50)
\put(-11,10){\psfig{figure=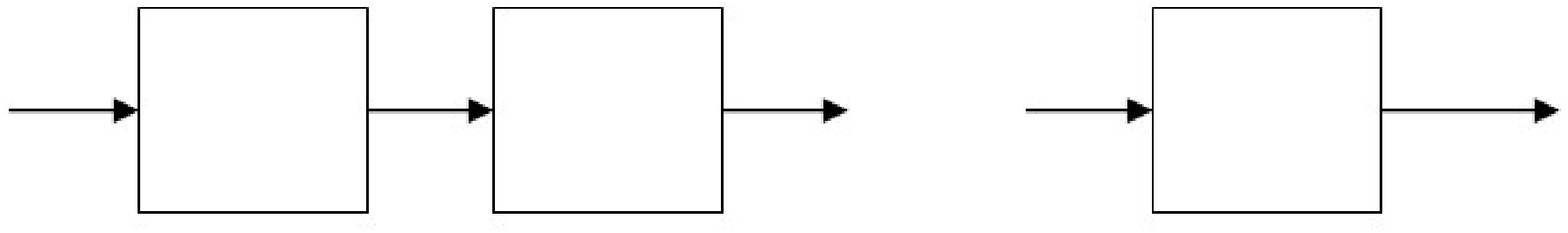,width=90\unitlength}}
\put(39,20){$\displaystyle \cdots$} \put(5,20){$\displaystyle H_1$}
\put(24,20){$\displaystyle H_2$} \put(58,20){$\displaystyle H_n$}
\put(-9,20){$\displaystyle u$}\put(77,20){$\displaystyle y_n$}
\end{picture}
\vspace{-1.5cm} \caption{The cascade interconnection for Corollary
\ref{two}.} \label{nest2}
\end{center}
\end{figure}

When the blocks $H_1,...,H_n$ each satisfy the OFP
property (\ref{ofp}), their cascade interconnection in Figure
\ref{nest2} inherits the sum of their phases and loses passivity.
The following corollary to Theorem \ref{main} quantifies the
``shortage" of passivity in such a cascade:
\begin{corollary}\label{two}
Consider the cascade interconnection in Figure \ref{nest2}. If each
block $H_i$ satisfies (\ref{ofp}) with a $C^1$ storage function
$V_i$ and $\gamma_i>0$, then for any
\begin{equation}\label{ifpcond}
\delta>\gamma_1\cdots \gamma_n \cos(\pi/(n+1))^{(n+1)},
\end{equation}
the cascade admits a storage function of the form (\ref{lyap})
satisfying
\begin{equation}\label{ifp}
\dot{V}\le -\epsilon |y|^2+\delta u^2 +uy_n.
\end{equation}
for some $\epsilon>0$. \hfill $\Box$
\end{corollary}

Inequality (\ref{ifp}) is an {\it input feedforward passivity} (IFP)
property \cite{sepulchre} where the number $\delta$ represents the
gain with which a feedforward path, if added from $u$ to $y_n$ in
Figure \ref{nest2}, would achieve passivity. Corollary \ref{two}
thus shows that the cascade of OFP systems (\ref{ofp}) in which
$\gamma_i>0$ represents an ``excess" of passivity, satisfies the IFP
property (\ref{ifp}) with a ``shortage"  characterized by
(\ref{ifpcond}).
\bigskip

\noindent {\bf Proof of Corollary \ref{two}:} Using (\ref{lyap}),
(\ref{ofp}), and substituting $u_i=y_{i-1}, \ i=2,\cdots,n,$  we
rewrite (\ref{ifp}) as
\begin{equation}\label{ineq}
d_1(-y_1^2+\gamma_1y_1u)+\sum_{i=2}^n
d_i(-y_i^2+\gamma_iy_iy_{i-1})+\delta(-u^2-\frac{1}{\delta} uy_n)\le
-\epsilon |y|^2.
\end{equation}
To show that $d_i>0,$ $i=1,\cdots,n,$ satisfying (\ref{ineq}) indeed
exist, we define
\begin{equation}\label{Atildematrix}
\tilde{A}=\left[\begin{array}{ccccc} -1  & 0 & \cdots & 0 & -\frac{1}{\delta} \\ \gamma_1 & -1 &\ddots & & 0 \\
0 & \gamma_2 & -1 & \ddots & \vdots \\ \vdots & \ddots & \ddots &
\ddots & 0 \\ 0 & \cdots & 0 & \gamma_n & -1 \end{array} \right]
\end{equation}
and note that the left-hand side of (\ref{ineq}) is
\begin{equation}\label{ineq2}
[u\ y^T]\tilde{D}\tilde{A}\left[\begin{array}{c}u
\\ y\end{array}\right]
\end{equation}
where $\tilde{D}:={\rm diag}\left\{\delta, \, d_1,\cdots, d_n
\right\}$. Because $\tilde{A}$ is of the form (\ref{Amatrix}) with
dimension $(n+1)$, an application of Theorem \ref{main} shows that a
diagonal $\tilde{D}$ rendering (\ref{ineq2}) negative definite
exists if and only if
$(\gamma_1\cdots\gamma_n\frac{1}{\delta})<\sec(\pi/(n+1))^{(n+1)}$.
Because this condition is satisfied when $\delta$ is as in
(\ref{ifpcond}), we conclude that such a $\tilde{D}>0$ exists and,
thus, (\ref{ifp}) holds.
 \hfill $\Box$

 \subsection*{\bf APPENDIX: Proof of Theorem \ref{main}}
Necessity follows because, as shown in \cite{TysOth78},  if
(\ref{eduardo}) fails then $A$ is not Hurwitz. To prove that
(\ref{eduardo}) is sufficient for diagonal stability, we define
\begin{eqnarray}
r&:=&(\gamma_1\dots\gamma_n)^{1/n}>0 \\ \Delta&:=&{\rm diag}
\left\{1,\ -\frac{\gamma_2}{r},\ \frac{\gamma_2\gamma_3}{r^2},\
\cdots, (-1)^{i+1}\frac{\gamma_2\cdots \gamma_i}{r^{i-1}}, \cdots,
(1)^{n+1}\frac{\gamma_2\cdots \gamma_n}{r^{n-1}} \right\} \nonumber
\end{eqnarray}
 and note that
\begin{equation}\label{structure}
-\Delta^{-1}A\Delta=\left[\begin{array}{c}\begin{array}{ccccc} 1  & 0 & \cdots & 0 & (-1)^{n+1}r \\ r & 1 &\ddots & & 0 \\
0 & r & 1 & \ddots & \vdots \\ \vdots & \ddots & \ddots & \ddots & 0
\\ 0 & \cdots & 0 & r & 1 \end{array}
\end{array}\right].
\end{equation}
Thus, with the choice
\begin{equation}
D=\Delta^{-2}
\end{equation}
we get
\begin{equation}\label{sym}DA+A^TD=\Delta^{-1}(\Delta^{-1}A\Delta+\Delta A^T
\Delta^{-1})\Delta^{-1} \end{equation} which means that $ DA+A^TD<0$
holds if the symmetric part of (\ref{structure}), given by
\begin{equation}\frac{1}{2}(-\Delta^{-1}A\Delta-\Delta A^T
\Delta^{-1}), \label{str2}
\end{equation}
is positive definite. To show that this is indeed the case, we note
that (\ref{structure}) exhibits a {\it circulant} structure
\cite{davis} when $n$ is odd, and a {\it skew-circulant} structure
when $n$ is even. In particular, it admits the
eigenvalue-eigenvector pairs
$$
\lambda_k=1+re^{i\frac{2\pi}{n}k} \quad v_k=\frac{1}{n}[1 \
e^{-i\frac{2\pi}{n}k} \ e^{-i2\frac{2\pi}{n}k} \cdots
e^{-i(n-1)\frac{2\pi}{n}k}]^T \quad k=1,\cdots,n
$$
when $n$ is odd; and
$$
\lambda_k=1+re^{i(\frac{\pi}{n}+\frac{2\pi}{n}k)} \quad
v_k=\frac{1}{n}[1 \ e^{-i(\frac{\pi}{n}+\frac{2\pi}{n}k)} \
e^{-i2(\frac{\pi}{n}+\frac{2\pi}{n}k)} \cdots
e^{-i(n-1)(\frac{\pi}{n}+\frac{2\pi}{n}k)}]^T
$$
when $n$ is even. Since, in either case, (\ref{structure}) is
diagonalizable with the unitary matrix $V=[v_1\cdots v_n]$, the
eigenvalues of the symmetric part (\ref{str2}) coincide with the
real parts of $\lambda_k$'s above. Finally, because
$$
\min_{k=1,\cdots n}{Re}\{1+re^{i\frac{2\pi}{n}k}\}=\min_{k=1,\cdots
n}{Re}\{1+re^{i(\frac{\pi}{n}+\frac{2\pi}{n}k)}\}=1-r\cos(\pi/n),
$$
we conclude that if (\ref{eduardo}) holds, that is $r<\sec(\pi/n)$,
then all eigenvalues of (\ref{str2}) are positive and, hence,
(\ref{str2}) is positive definite and (\ref{sym}) is negative
definite. \hfill $\Box$

\end{document}